\documentclass[12pt]{amsart}
\usepackage{amsmath,amssymb,amsfonts,amsthm}
\usepackage{epsfig}
\usepackage{a4wide}
\usepackage{enumerate,color}
\usepackage{subfigure,hyperref}
\usepackage{aliascnt}
\usepackage{todonotes}

\newtheorem{theorem}{Theorem}

\newaliascnt{lemma}{theorem}
\newtheorem{lemma}[lemma]{Lemma}
\aliascntresetthe{lemma}

\newaliascnt{prob}{theorem}
\newtheorem{prob}[prob]{Problem}
\aliascntresetthe{prob}

\newaliascnt{conjecture}{theorem}
\newtheorem{conjecture}[conjecture]{Conjecture}
\aliascntresetthe{conjecture}

\newaliascnt{observation}{theorem}
\newtheorem{observation}[observation]{Observation}
\aliascntresetthe{observation}

\hypersetup{
  pdftitle = {Polynomial expansion and sublinear separators},
  pdfauthor = {Louis Esperet and Jean-Florent Raymond}
}

\DeclareMathOperator{\tw}{tw}
\DeclareMathOperator{\polylog}{polylog}

\def\cC{\mathcal{C}}

\def\longequation{$$\vcenter\bgroup\advance\hsize by -9em%
\noindent\ignorespaces\refstepcounter{equation}}%
\makeatletter%
\def\endlongequation{\egroup\eqno(\theequation)$$\global\@ignoretrue}
\makeatother

\let\oldrceil\rceil
\renewcommand{\rceil}{\right\oldrceil}
\let\oldlceil\lceil
\renewcommand{\lceil}{\left\oldlceil}

\begin{document}

\title{Polynomial expansion and sublinear separators}

\author[L. Esperet]{Louis Esperet} 
\address{Univ. Grenoble Alpes, CNRS, G-SCOP, Grenoble, France}
\email{louis.esperet@grenoble-inp.fr}

\author[J.-F. Raymond]{Jean-Florent Raymond}
\address{Institute of Informatics, University of Warsaw, Poland and LIRMM, University of Montpellier, France}
\email{jean-florent.raymond@mimuw.edu.pl}

\thanks{The first author was supported by ANR Projects STINT
  (\textsc{anr-13-bs02-0007}) and GATO (\textsc{anr-16-ce40-0009-01}) and LabEx PERSYVAL-Lab
  (\textsc{anr-11-labx-0025-01}). The second author was supported by the Polish National Science Centre grant PRELUDIUM DEC-2013/11/N/ST6/02706.}

\date{}

\sloppy

\begin{abstract}
Let $\cC$ be a class of graphs that is closed under taking
subgraphs. We prove that if for some fixed
$0<\delta\le 1$, every $n$-vertex graph of $\cC$ has
a balanced separator of order $O(n^{1-\delta})$, then any depth-$r$ minor (i.e.\ minor obtained by
contracting disjoint subgraphs of radius at most $r$) of
a graph in $\cC$ has average degree $O\big((r \polylog 
r)^{1/\delta}\big)$. This confirms a conjecture of
Dvo\v r\'ak and Norin.
\end{abstract}

\maketitle

\section{Introduction}\label{sec:intro}

For an integer $r\ge 0$, a \emph{depth-$r$ minor} of a graph $G$ is a
subgraph of a graph that can be obtained from $G$ by contracting pairwise
vertex-disjoint subgraphs of radius at most $r$.
Let $d(G)$ denote the average degree of a graph $G=(V,E)$, i.e.\
$d(G)=2|E|/|V|$. 
For some function $f$, we say that a class $\cC$ of graphs has \emph{expansion bounded
by $f$} if for any graph $G\in \cC$ and any integer $r$, any depth-$r$
minor of $G$ has average degree at most $f(r)$. We say that a class
has \emph{bounded expansion} if it has expansion bounded by some function
$f$, and \emph{polynomial expansion} if $f$ can be taken to be a
polynomial.

\smallskip

Classes of bounded expansion play a central role in the study of
sparse graphs~\cite{NO12}. From an algorithmic point of view, a very
useful property of theses classes is that when their expansion is not too
large (say subexponential), graphs in the class have sublinear
separators.
A \emph{separator} in a graph $G=(V,E)$ is a pair of subsets $(A,B)$ of
vertices of $G$ such that $A\cup B =V$ and no edge of $G$ has one endpoint in
$A\setminus B$ and the other in $B \setminus A$. The separator $(A,B)$
is said to be \emph{balanced} if both $|A\setminus B|$ and $|B\setminus A|$ contain at most
$\tfrac23 |V|$ vertices. The \emph{order} of the
separator $(A,B)$ is $|A\cap B|$. 

\smallskip

A class $\cC$ of graphs is \emph{monotone} if for any graph $G \in \cC$, any
subgraph of $G$ is in $\cC$. Dvo\v r\'ak and Norin~\cite{DN16} observed that the following can be
deduced from a result of Plotkin, Rao, and Smith~\cite{PRS94}.

\begin{theorem}[\cite{DN16}]\label{thm:prs}
Let $\cC$ be a monotone class of graphs with expansion bounded by
$r\mapsto c (r+1)^{1/4\delta-1}$, for some constant $c>0$ and $0<\delta\le 1$. Then there is
a constant $C$ such that every
$n$-vertex graph of $\cC$ has a balanced separator of order $Cn^{1-\delta}$.
\end{theorem}

Dvo\v r\'ak and Norin~\cite{DN16} also proved the following partial converse.

\begin{theorem}[\cite{DN16}]\label{thm:zs}
Let $\cC$ be a monotone class of graphs such that for some fixed constants $C>0$ and $0<\delta\le 1$, every
$n$-vertex graph of $\cC$ has a balanced separator of order $C
n^{1-\delta}$. Then the expansion of $\cC$ is
bounded by some function $f(r) = O(r^{5/\delta^2})$.
\end{theorem}

They conjectured that the exponent $5/\delta^2$ of the polynomial expansion
in \autoref{thm:zs} could be improved to match (asymptotically) that of
\autoref{thm:prs}. 

\begin{conjecture}[\cite{DN16}]\label{conj:zs}
There exists a real $c>0$ such that the following holds. Let $\cC$ be a monotone class of graphs such that for some fixed constants $C>0$ and $0<\delta\le 1$, every
$n$-vertex graph of $\cC$ has a balanced separator of order $C
n^{1-\delta}$. Then the expansion of $\cC$ is
bounded by some function $f(r) = O(r^{c/\delta})$.
\end{conjecture}

In this short note, we prove this conjecture.

\begin{theorem}\label{thm:main}
For any $C>0$ and
$0<\delta\le 1$, if a monotone class $\cC$ has the property that every
$n$-vertex graph in $\cC$ has a balanced separator of order at most $C
n^{1-\delta}$, then $\cC$ has expansion bounded by the function $f:
r \mapsto c_1\cdot (r+1) ^{1/\delta} (\tfrac1{\delta}\log
(r+3))^{c_2/\delta}$, for some constants $c_1$ and $c_2$
depending only on $\cC$.
\end{theorem}

In particular \autoref{conj:zs} holds for any real number $c> 1$. The proof of \autoref{thm:main} is given in the next section, and we
conclude with some open problems in \autoref{sec:pro}.

\section{Proof of \autoref{thm:main}}\label{sec:proof}

We need the following results. The first is a classical
connection between balanced separators and tree-width (see~\cite{DN16}).

\begin{lemma}\label{lem:treewidth}
  Any graph $G$ has a balanced separator of order at most $\tw(G)+1$.
\end{lemma}

Dvo\v r\'ak and Norin~\cite{DN14} proved that the following partial converse holds.

\begin{theorem}[\cite{DN14}]\label{thm:treewidth}
If every subgraph of $G$ has a balanced separator of order at most $k$,
then $G$ has tree-width at most $105k$.
\end{theorem}

Note that in our proof of \autoref{thm:main} we could also use the
weaker (and easier) result of~\cite{BGHK95} that under the same
hypothesis, $G$ has tree-width at most $1+k \log |V(G)|$, but the
computation is somewhat less cumbersome if we use
\autoref{thm:treewidth} instead.

\medskip

For a set $S$ of vertices in a graph $G$, we let $N(S)$ denote the set of
vertices not in $S$ with at least one neighbor in $S$. We will use the
following result of Shapira and Sudakov~\cite{SS15}.

\begin{theorem}[\cite{SS15}]\label{thm:expander}
Any graph $G$ contains a subgraph $H$ of average degree $d(H)\ge
\tfrac{255}{256} d(G)$ such that for any set $S$
of at most $n/2$ vertices of $H$ (where $n=|V(H)|$), $|N(S)|\ge
\tfrac1{2^8 \log n (\log \log n)^2}|S|$.
\end{theorem}

In fact, we will only need a much weaker version, where the
vertex-expansion is of order $\Omega \left (\tfrac1{\polylog  n}\right )$ instead of
$\Omega\left (\tfrac1{\log n(\log \log n)^2} \right)$.

\medskip

Finally, we need a result of Chekuri and Chuzhoy~\cite{CC15} on bounded-degree subgraphs of large
tree-width in a graph of large tree-width.

\begin{theorem}[\cite{CC15}]\label{thm:sparsifier}
There are constants $\alpha,\beta$ such that for any integer $k\geq 2$, any graph $G$ of tree-width at least $k$ contains
a subgraph $H$ of tree-width at least $\alpha k/(\log k)^\beta$ and
maximum degree~3.
\end{theorem}

Let us remark that instead of \autoref{thm:sparsifier}, our proof of
\autoref{thm:main} could rely on an earlier result of Chekuri and
Chuzhoy~\cite{ChekuriC13larg} which, under the same assumptions,
merely guarantees the existence of a subgraph of $G$ of treewidth
$\Omega(k/(\log k)^6)$ and maximum degree $O((\log k)^3)$.

\medskip

We are now ready to prove our main result.

\medskip

\noindent \emph{Proof of \autoref{thm:main}.} Let $G$ be a graph of
$\cC$ and let $F$ be a depth-$r$ minor of $G$. Our goal is to prove
that $d(F)\le c_1\cdot (r+1) ^{1/\delta} (\tfrac1{\delta}\log
(r+3))^{c_2/\delta}$, for some constants $c_1$ and $c_2$ depending
only on $C$. Note that for any $r\ge 0$ and $0<\delta\le
1$,
\[
  c_1\cdot (r+1) ^{1/\delta} (\tfrac1{\delta}\log (r+3))^{c_2/\delta}\ge \max\left\{c_1(\log 3)^{c_2},c_1
    \exp(\tfrac{c_2}{\delta}\log\tfrac{\log 3}{\delta})\right\},
\]
so we can assume without loss of generality that
\[
  d(F)\ge \max\left\{{10}^8,
    \exp\left(4\cdot\tfrac{\beta+3}{\delta}\log(2\cdot\tfrac{\beta+3}{\delta})\right)\right\}
\]
by choosing appropriate values of $c_1,c_2$. By
\autoref{thm:expander}, $F$ has a subgraph $H$ of average degree
$d(H)\ge \tfrac{255}{256} d(F)$ such that for any set $S$ of at most
$|V(H)|/2$ vertices of $H$,
\[
  |N(S)|\ge \tfrac1{2^8 \log |V(H)|  (\log \log |V(H)| )^2}|S|\ge
  \tfrac1{2^8 (\log |V(H)|)^3}|S|.
\]

 It follows from \autoref{lem:treewidth} that $H$ contains a balanced
 separator $(A,B)$ with $|A\cap B|\le\tw(H)+1$. As $A\setminus B$ and
 $B\setminus A$ are disjoint, one of them contains at most half of the
 vertices. We may assume without loss of generality that $|A\setminus
 B|\le |V(H)|/2$. As $N(A \setminus B) \subseteq A \cap B$, we get
\[
  |A\cap B|\ge \tfrac1{2^8 (\log
    |V(H)|)^3  }|A\setminus B|.
\]
Since $(A,B)$ is
balanced, $|A\setminus B|+|A\cap B|\ge \tfrac13 |V(H)|$ and so
\[
  \tfrac13|V(H)|\le |A\cap B|(1+2^8( \log |V(H)|)^3).
\]
Given that $|A\cap B|\le\tw(H)+1$, we deduce
\[
  \tw(H)\ge \tfrac{|V(H)|}{3\cdot 2^{8} (\log  |V(H)|)^3+3}-1 \ge
  \tfrac{|V(H)|}{2^{10} (\log |V(H)|)^3},
\]
using that $|V(H)|\ge d(H) \ge \tfrac{255}{256} \cdot 10^8$.

\medskip

By \autoref{thm:sparsifier}, $H$ has a subgraph $H'$
of maximum degree 3 such that
\[
  \tw(H')\ge \tfrac{\alpha \,\tw(H)}{(\log
\tw(H))^\beta}\ge \tfrac{\alpha \,|V(H)|}{2^{10} (\log
|V(H)|)^{\beta+3} },
\]
since $\tw(H)\le |V(H)|$. Note that $H'$ is a subgraph
of $H$ (and $F$) and therefore also a depth-$r$ minor of $G$. In $G$, $H'$ corresponds to a
subgraph $G'$ (before contraction of the subgraphs of radius $r$) with
$|V(G')|\le (3r+1)|V(H')|\le (3r+1)|V(H)|$. Indeed, since $H'$ has maximum degree 3, each
subgraph of radius at most $r$ in $G'$ whose contraction corresponds to a
vertex of $H'$ contains at most $3r+1$ vertices. Since $H'$ is a
minor of $G'$, we have
\[
  \tw(G')\ge \tw(H')\ge \tfrac{\alpha \,|V(H)|}{2^{10} (\log
    |V(H)|)^{\beta+3} }.
\]

Since $\cC$ is monotone, every subgraph of $G'$ is in $\cC$ and thus has a balanced
separator of order at most $C |V(G')|^{1-\delta}$. Hence, by
\autoref{thm:treewidth},
\[
  \tw(G')\le 105 C |V(G')|^{1-\delta}\le 2^7 C |V(G')|^{1-\delta}.\label{eq:upperbound}
\]

We just obtained lower and upper bounds on
$\tw(G')$. Putting them together, we obtain:

\begin{align*}
  \tfrac{\alpha \,|V(H)|}{2^{10} (\log
    |V(H)|)^{\beta+3} } &\le 2^7 C\, 
                          |V(G')|^{1-\delta}\\
                        &\le 2^7 C\,
  \big((3r+1)|V(H)|\big)^{1-\delta}\text{, and thus}\\
  \tfrac{|V(H)|^\delta}{(\log |V(H)|)^{\beta+3}} &\le \tfrac{2^{17}
                                                   C}{\alpha}(3r+1)^{1-\delta}\\
                        &\le \tfrac{2^{17}C}{\alpha}(3r+1)
  \text{, and}\\
  |V(H)| &\le \left( \tfrac{2^{17}C}{\alpha}(3r+1) (\log
    |V(H)|)^{\beta+3}\right)^{1/\delta}.
\end{align*}

It follows that
\[
  \log|V(H)|\le \tfrac1{\delta} \log
  \left (\tfrac{2^{17}C}{\alpha}(3r+1)\right) +\tfrac{\beta+3}{\delta} \log
  \log |V(H)|.
\]

Since the function $n\mapsto \tfrac{\log n}{\log \log n}$ is increasing
for $n\ge 16$, a direct consequence of our initial
assumption that $|V(H)|\ge
\exp\left(4\cdot\tfrac{\beta+3}{\delta}\log(2\cdot\tfrac{\beta+3}{\delta})\right)$ is that

\begin{align*}
  \tfrac{\log |V(H)|}{\log \log |V(H)|} &\ge
                                          2\cdot\tfrac{\beta+3}{\delta}\text{,
                                          and thus}\\
  \log |V(H)|&\le \tfrac2{\delta} \log \left ( \tfrac{2^{17}C}{\alpha}(3r+1)\right ).                                        
\end{align*}

We conclude that $$|V(H)|\le \left( \tfrac{2^{17}C}{\alpha}(3r+1)
  \left( \tfrac2{\delta} \log \left (
\tfrac{2^{17}C}{\alpha}(3r+1)\right)\right)^{\beta+3}
\right)^{1/\delta}\le \tfrac{255}{256} c_1 (r+1)^{1/\delta} (\tfrac1{\delta}\log (r+3))^{c_2/\delta},  $$ for some constants $c_1,c_2$ depending only on
$C$ and the constants $\alpha,\beta$ of
\autoref{thm:sparsifier}. Recall that $d(F)\le
\tfrac{256}{255}d(H)$. Since $d(H)\le |V(H)|$, we obtain $d(F)\le c_1\cdot (r+1) ^{1/\delta} (\tfrac1{\delta}\log
(r+3))^{c_2/\delta}$, as desired. This concludes the proof of
\autoref{thm:main}. \hfill $\Box$

\section{Open problems}\label{sec:pro}

A natural problem is to determine the infimum real $c>0$, such that if a monotone class $\cC$ has the property that every
$n$-vertex graph in $\cC$ has a balanced separator of order
$O(n^{1-\delta})$, then $\cC$ has expansion bounded by some function
$r \mapsto O(r^{c/\delta})$. \autoref{thm:main}
implies that $c \le 1$. On the other hand, it directly follows from
\autoref{thm:prs} that $c\le \tfrac1{4+\epsilon}$ would imply that if any
$n$-vertex graph in $\cC$ has a balanced separator of order
$O(n^{1-\delta})$, then any
$n$-vertex graph in $\cC$ has a balanced separator of order
$O(n^{1-(1+\epsilon/4)\delta})$. Therefore, \autoref{thm:prs}
implies that $c \ge \tfrac14$ (moreover, the
proof of \autoref{thm:prs} in~\cite{DN16} can be slightly optimized to
show that $c \ge \tfrac12$). A good candidate to prove a better lower
bound for $c$ would be the family of all finite subgraphs of the
infinite $d$-dimensional
grid. The $n$-vertex graphs in this class have balanced separators of order
$O(n^{1-1/d})$ (see~\cite{MTV91}), and it might be the case
that they have expansion $\Omega(r^{cd})$ for some $c>\tfrac12$.

\medskip

One way to measure the sparsity of a class of graphs is via its expansion
(as defined in \autoref{sec:intro}). Another way (which turns out to
be equivalent) is via its \emph{generalized coloring
  parameters}. Given a linear order $L$ on the vertices of a graph
$G$, and an integer $r$, we say that a vertex $v$ of $G$ is \emph{strongly $r$-reachable} from a
vertex $u$ (with respect to $L$) if $v \le_L u$, and there is a path $P$ of length at most $r$  between $u$ and
$v$, such that $u <_L w$ for any internal vertex $w$ of $P$. If we
only require that $v$ is the minimum of the vertices of $P$ (with
respect to $L$), we say that $v$ is \emph{weakly $r$-reachable} from
$u$. The \emph{strong $r$-coloring number} $\text{col}_r(G)$ of $G$
is the minimum integer $k$ such that there is a linear order $L$ on the
vertices of $G$ such that for any vertex $u$ of $G$, at most $k$
vertices are strongly $r$-reachable from $u$ (with respect to $L$). By
replacing \emph{strongly} by \emph{weakly} in the previous definition,
we obtain the \emph{weak $r$-coloring number} $\text{wcol}_r(G)$ of
$G$. Note that for any graph $G$ and any integer $r$,
$\text{col}_r(G)\le \text{wcol}_r(G)$. For more on these parameters and their connections with the
expansion of graph classes,
the reader is referred to~\cite{NO12}.

As we have seen before, it follows from~\cite{DN16} that a monotone
class of graphs has polynomial expansion if and only if, for some fixed $0< \delta \le
1$,
each $n$-vertex graph in the class has a balanced
separator of order $O(n^{1-\delta})$. Joret and Wood asked whether this is also equivalent to having
weak and strong $r$-coloring numbers bounded by a polynomial function
of $r$.

\begin{prob}[Joret and Wood, 2017]\label{pro}
Assume that $\cC$ is a monotone class of graphs. Are the following statements
equivalent?
\begin{enumerate}\item $\cC$ has polynomial expansion.
\item There exists a constant $c$, such that for every $r$, every graph in $\cC$ has strong $r$-coloring number at most $O(r^c)$. 
\item There exists a constant $c$, such that for every $r$, every
  graph in $\cC$ has weak $r$-coloring number at most $O(r^c)$. 
\end{enumerate}
\end{prob}

Note that clearly (3) implies (2). It was known that (3) implies (1)
(this is a consequence of Lemma 7.11 in~\cite{NO12}), and Norin
recently made the following observation, which shows
that (2) implies (1). 

\begin{observation}[Norin, 2017]\label{obs:sergey} 
Every depth-$r$ minor of a graph $G$ has average degree at most
$2\,\mbox{\rm{col}}_{4r}(G)$. 
\end{observation}

\begin{proof}
Let $L$ be a linear order on the vertices of $G$, such that for any vertex $v$ of $G$, at most $\text{col}_{r}(G)$ vertices are strongly $r$-reachable from $v$ (with respect to $L$). Let $H$ be a depth-$r$ minor of a graph $G$. For any vertex $u$ of $H$, let $S_u$ be a
subgraph of $G$ of radius at most $r$, such that the $S_u$'s are
vertex-disjoint and for any edge $uv$ of $H$, there is an edge in $G$
between a vertex of $S_u$ and a vertex of $S_v$.  It is enough to prove that there is a linear order $L'$ on the vertices
of $H$ such that any vertex $u$ of $H$, at most $\text{col}_{4r}(G)$
vertices of $H$ are strongly 1-reachable from $u$.

We construct $L'$
from $L$ as follows: for $u,v$ in $H$, we set $u<_{L'} v$ if and only if,
with respect to $L$, the smallest vertex of $S_u$ precedes the
smallest vertex of $S_v$. This clearly defines a linear order on the
vertices of $H$. Consider a vertex $u$ of $H$ and let $x$ be the smallest vertex of
$S_u$ (with respect to $L$). Let $v$ be a neighbor of
$u$ in $H$ with $v<_{L'} u$ (i.e.\ $v$ is strongly 1-reachable from
$u$ in $H$). Let $t\in S_u$ and $z\in S_v$ be such that $tz$ is an edge of
$G$. Observe that there is a path $P_u$ from $x$ to $t$ in $S_u$ (and
$x$ is the smallest vertex in this path with respect to $L$), and a
path $P_v$ from $z$ to $y$ in $S_v$. Let $w$ be the first vertex of
$P_v$ such that $w<_L x$ (note that possibly $w=z$). The concatenation
of $P_u$, $zt$, and the subpath of $P_v$ between $z$ and $w$ has length at most $4r$ and thus shows that $w$
is strongly $4r$-reachable from $x$ in $G$. Hence, at most $\text{col}_{4r}(G)$
vertices of $H$ are strongly 1-reachable from $u$ in $H$ with respect
to $L'$, as desired. 
\end{proof}

\section*{Acknowledgements}

We thank Zden\v ek Dvo\v r\'ak for the discussion about~\cite{DN16}, 
Gwena\"el Joret and David Wood for allowing us to mention \autoref{pro}, and
Sergey Norin for allowing us to mention \autoref{obs:sergey} and
its proof.

\end{document}